\documentclass{amsart}
\usepackage[cp1251]{inputenc}
\usepackage[english,russian]{babel}
\usepackage{amsmath}
\usepackage{amssymb}
\usepackage{amsfonts}
\def\udcs{517.98} 

\newtheorem{lemma}{Lemma}
\newtheorem{theorem}{Theorem}
\newtheorem{definition}{Definition}

\def\Re{\operatorname{\mathrm Re}}
\def\Im{\operatorname{\mathrm Im}}
\def\Fin{\operatorname{\mathrm Fin}}
\def\Ker{\operatorname{\mathrm Ker}}

\def\exp{\operatorname{\mathrm exp}}
\def\integer{\operatorname{\mathrm Int}}
\def\Fin{\operatorname{\mathrm Fin}}
\def\bR{\mathbb R}
\def\bS{\mathbb S}
\def\bC{\mathbb C}
\def\bN{\mathbb N}
\def\bZ{\mathbb Z}
\begin{document}
УДК \udcs
\thispagestyle{empty}

\title[Interpolation by sums of the series of exponentials in $H (\mathbb C)$\dots ]
{Interpolation by sums of the series of exponentials in $H (\mathbb C)$ with interpolation nodes on the rays.} 
\author{S.G. Merzlyakov, S.V. Popenov}
\address{Sergey Victorovich Popenov, 
\newline\hphantom{iii} Institute of Mathematics, Ufa Scientific Center, Russian Academy of Sciences,
\newline\hphantom{iii} Chernishevskii str., 112,  
\newline\hphantom{iii} Russian Federation, 450008, Ufa}  
\email{msg2000@mail.ru}
\address{Sergey Georgievich Merzlyakov, 
\newline\hphantom{iii} Institute of Mathematics, Ufa Scientific Center, Russian Academy of Sciences,
\newline\hphantom{iii} Chernishevskii str., 112,  
\newline\hphantom{iii} Russian Federation, 450008, Ufa}  
\email{spopenov@gmail.com}

\thanks{\rm This work was supported by the Russian Foundation for Basic Research [grant number 15-01-01661 А]; }
\maketitle 
{
\small
\begin{quote}

\noindent{\bf Abstract. } In the space of all entire functions it is solved the problem of interpolation taking into account multiplicities by sums of the series of  exponentials with the exponents from a given set.  It is found a criterion of solubility of the problem in the case when all infinite portions of interpolation nodes are situated on a finite system of rays. It is also disclosed that the problem is equivalent to particular problems of simple interpolation, as well as of point-wise crude approximation, by sums of the series of exponentials.
 
\end{quote}
}

\section{The series of exponentials with the exponents from a given set}
Let $H(\mathbb C)$ be the space of all entire functions. If $\{K_j\},\,j\in \bN$, is a sequence of convex compact sets such that $\bC=\cup K_j,\,K_j\subset \integer K_{j+1}$, the topology in $H(\bC)$ is defined by the norms 
  \begin{equation*}
  p_{K_j}(f)= \max_{t\in K_j}|f(t)|,\,f\in H(\bC).
  \end{equation*}

Let  $\Lambda$ be some fixed infinite subset of the complex plane. In order not to burden the demonstration, suppose herein that the set of all finite limit points of $\Lambda$ is bounded. Denote
 \begin{equation*} \Sigma(\bC,\Lambda)=\{u: u(z)=\sum_{n=1}^{\infty} c_{n} e^{\lambda_nz},\,z\in \bC,\,\{\lambda_n\}\subset \Lambda\},
 \end{equation*}
wherein it is supposed that the series of exponentials is absolutely convergent for all $z \in\bC$. 

 The set $\Sigma(\bC,\Lambda)$ is not empty. Indeed, for arbitrary sequence $\{\lambda_j\}$ and for every $j\in \bN$ take $c_j$ so as to $|c_j|p_{K_j}(e^{\lambda_j\cdot})\leqslant 1/2^j$, then for all $z\in K_j$  and $n\geqslant j+1$, $|c_n|p_{K_j}(e^{\lambda_n\cdot})\leqslant 1/2^n$. Thus the series $\sum_{n=j+1}^{\infty}c_n e^{\lambda_nz}$ is absolutely convergent as well as the series of norms. In particular such a series is convergent in the topology of the space $H(D)$. 
 
 It is known  (see, e.g., \cite{Leonu}, \cite{Mcont}) that if series of exponentials converges absolutely in a convex domain it converges uniformly in it. Moreover the series of exponentials is normally convergent that is for any convex compact set $K\subset \bC$ the numerical series of norms
\begin{equation*}
\sum_{n=1}^{\infty}|c_n| p_K(e^{\lambda_n\cdot})=\sum_{n=1}^{\infty}|c_n|e^ {H_K(\lambda_n)}
\end{equation*}
 is convergent, wherein $H_K(z)=\sup_{\sigma\in K} \Re (\sigma z),\, z\in \bC$, is the support function of the compact set $K$ in the sense of complex analysis. Hence every series in the definition of $\Sigma(\bC,\Lambda)$ is uniformly convergent on all compact subsets of $\bC$ and $ \Sigma(\bC,\Lambda)$ is a linear subspace of $H(\bC)$. In general it is not closed subspace.

It may be noted that $H_K(z)=h_K(\theta)|z|,\,z=|z|e^{i\theta}$, where $h_K(\theta)=\sup_{\sigma\in K} \Re (\sigma e^{i\theta})$ is the support function (in the sense of $\bR^2$) of the compact set $\overline K$ that is complex conjugate of $K$. 

\medskip
\noindent{\bf Problem of interpolation by sums of the series of exponentials in $H (\mathbb C)$.}
 Denote by $\mathcal M =\{\mu_k\in \bC,\, k\in \mathbb N\}$ the set of interpolation nodes taking into account multiplicities $m_k,\,k\in \mathbb N$. It is supposed that the set is discrete in $\bC$ i.e. it don`t have any finite limit points. 
 
  For a given $\Lambda$ in the space $H(\bC)$ consider the following problem of interpolation taking into account multiplicities, by sums of series of exponentials with the exponents from $\Lambda$:

{\noindent \it To describe discrete node sets $\mathcal M\subset \bC$ that allow, for arbitrary interpolation data $b^j_k\in \bC$, to find such a function $u\in \Sigma(\bC,\Lambda)$ that $u^{(j)}(\mu_k)=b^j_k,\,k\in \bN,\,j=0,\cdots,m_k-1$.} 

To be more precise, to address the problem one needs to describe both $\mathcal M$ and $\Lambda$. 

\medskip
{\bf Example.} Consider two nodes $\mu_1=0,\,\mu_2=\exp(i\varphi),\,0\leqslant\varphi<2\pi,\,$ and $\Lambda=\{\lambda_n=2\pi n \exp(i(\pi/2-\varphi)),\,n\in\bZ\}$. Then all the functions $u\in\Sigma(\bC,\Lambda)$ satisfy $f(\mu_1)=f(\mu_2)$ because of $\exp(\mu_2\lambda_n)=1$ and interpolation fails. The example provides a necessary condition for considered interpolation problem.

Further, for solubility of the problem it is also necessary to demand the set $\Lambda$  to be unbounded, otherwise any function $u\in \Sigma(\bC,\Lambda)$ has the continuation to the entire function of exponential type. 
Indeed, when  $|\lambda_n|\leqslant C$ and the series $u(z)=\sum_{n=1}^{\infty} c_{n} e^{\lambda_nz}$ converges absolutely in a point $z_0\in \bC$, the simple estimate $|e^{\lambda_nz}|\leqslant e^{C|z-z_0|} |e^{\lambda_nz_0}|,\,z \in\bC$, implies that $$|u(z)|\leqslant e^{C|z-z_0|}\sum_{n=1}^{\infty} |c_n||e^{\lambda_nz_0}|$$.
The estimate shows that simple interpolation is impossible for interpolation data having rapid growth (e.g., for $|b_k|\geqslant k e^{C|\mu_k-z_0|},\,z_0\in \bC$).

Conditions of solubility of the considered interpolation problem are formulated in the terms of limit directions at the infinity for the set  $\Lambda$.
 
Denote $\bS=\{s\in\bC: |s|=1\}$ and denote by $s=s_\alpha$ elements $s=e^{i\alpha}\in \bS$. 

\begin{definition}  Let $\Lambda\subset \bC$ be an arbitrary unbounded set. The set $P(\Lambda)\subset \bS$ of limit directions at infinity consists of all $s_\alpha\in \bS$ such that for some unbounded subsequence $\{\lambda_{n_k}\}\subset \Lambda,\ $ $s_\alpha=\lim_{|\lambda_{n_k}|\to\infty}\lambda_{n_k}/|\lambda_{n_k}|$. 
\end{definition}
Some simple properties: $ P(\Lambda)$ is the closed subset of $\bC$ and
\begin{equation*} 
P(\Lambda+\lambda_0)=P(\Lambda),\, P(rs_{\alpha}\Lambda)=s_{\alpha}P(\Lambda),\,\lambda_0\in\bC,\, r>0.
\end{equation*}

When $\mathcal M\subset \bR$ (or it is subset of some ray or direct line) in the space $H(\bC)$ the problem of interpolation by sums of the series of exponentials was completely solved (\cite{MP}).  Conditions were found on $\Lambda$ in the terms of location of some of limit directions from $P(\Lambda)$ that gives criteria of solubility of the problem. 

Denote by $I(\bC,\mathcal M)$ the closed ideal in $H(\bC)$ consisting of all entire functions $h$ that vanish on the set $\mathcal M$ with multiplicities $m_k$:
\begin{equation*}I(\bC,\mathcal M)=\{u^{(j)}(\mu_k)=0,\,k\in \bN,\,j=0,\cdots,m_k-1\}.
 \end{equation*}
 The classical theorems of complex analysis imply that the problem of interpolation by sums of the series of exponentials is equivalent to the existence of representations 
\begin{equation} \label{F1}
H(\bC)=\Sigma (\bC,\Lambda)+I(\bC,\mathcal M)
\end{equation}
corresponding to all discrete $\mathcal M\subset \bC$. Specifically this is a direct consequence of solubility of holomorphic interpolation problem \cite{Her} due to Weierstrass and Mittag-Leffler theorems.

\medskip
\noindent {\bf The problem of interpolation by elements of the kernel of convolution operator.}
In the paper \cite{Na} it was studied the multy-point Vall\' ee-Poussin problem for convolution operator $M_\varphi$, that is correctly defined in $H(\bC)$, with the nodes on the real line. Here $\varphi$ is an entire function of exponential type by which the operator $M_\varphi$ is generated. 

The problem may be considered as being the problem of interpolation by elements of the kernel of the convolution operator. More particularly V.V. Napalkov was the first who find sufficient conditions (\cite{NapP}, \cite{Na1},\cite{Na})  for existence of the Fischer type representations in the terms of zero set of function $\varphi$,
\begin{equation} \label{F2}
H(\bC)=\Ker M_\varphi+I(\bC,\mathcal M),
\end{equation} 
what is equivalent to the solubility of the said interpolation problem. We disclose some examples showing that the conditions of V.V. Napalkov are not the necessary ones for the problem. It is still the open problem of complete description of convolution operators in the setting.

The cited results may be derived from the mentioned criteria \cite{MP}, \cite{MP1} of interpolation by sums of series of exponentials while avoiding superfluous limitations on the function $\varphi$ prescribed in \cite{Na},  \cite{NZ}. In fact, lets denote by $Z_\varphi$ the zero set of the function $\varphi$  then all exponentials $\exp \lambda_n z\in \Ker M_\varphi$ for $\lambda_n\in Z_\varphi$. Define  $\Lambda$ as a subset of $Z_\varphi$. Then $\Sigma (\bC,\Lambda)\subset \Ker M_\varphi$ and if representation (\ref{F1}) is in effect for some node set $\mathcal M$ it implies the representation (\ref{F2}).

It is known (\cite{MP}, \cite{MP1}) that the representation of functions in (\ref{F1}) and (\ref{F2}) is not unique: if the problem is solvable the subspace $\Sigma (\bC,\Lambda) \cap I(\bC,\mathcal M)$ has infinite dimension. 

There are also results (see e.g., \cite{NZ}, \cite{MP1} and \cite{NN}) for the spaces $H(D)$ as well as for some complex sets of interpolation nodes, where $D\subset \bC$  is a convex domain. In \cite{MP1} the problem is solved  for $H(D)$ in arbitrary convex domain $D$ and for all node sets that are situated on the portion of some straight line in $D$ and have only finite limit points on the boundary. We are about to treat the situations in the short run. 
\section{The main results}
In the paper the said problem of interpolation, i.e. existence of the representations (\ref{F1}), is thoroughly settled for  the case when all infinite portions of the node sets $\mathcal M$ are situated on a system of rays $L_{\beta j}=\{z\in\bC: z=a_j+t\exp(i\beta_j),\, t>0\},\,$ $1\leqslant j\leqslant q,\,0\leqslant\beta_j<2\pi$. To be more precise, all the sets $ \mathcal M_j=\mathcal M\cap L_{\beta j}$ are infinite and there are at most finite number of $\mu_k\in\mathcal M$ outside $\cup_{j=1}^q\mathcal M_j$.
\begin{theorem}\label{T1}
Let $\Lambda$ be a given unbounded set. The problem of interpolation by functions from $\Sigma(\Lambda, \bC)$ with arbitrary discrete node set $\mathcal M$ being situated on the system of the rays as described above, is solvable  if and only if $\mathcal M$ and $\Lambda$ possess the properties:

(i) There exists a plurality $s_{\alpha_{\nu(j)}}\in P(\Lambda)$, say $1\leqslant \nu(j)\leqslant \tau\leqslant q$, such that   $|\beta_j+\alpha_\nu(j)|<\pi/2$ for every $j,\,j=1,\cdots,q;$ and

(ii)  	For every $j,\,1\leqslant j\leqslant q,\,$ and every $\mu_l\in\mathcal M_j,\,$  $\Re\mu_k\exp i\alpha_{\nu(j)}\not=\Re\mu_l\exp i\alpha_{\nu(j)}$ for all $k\in\bN,\, k\not=l$.
\end{theorem} 
It may be noted that (i) and (ii) describes the following transparent geometrical configurations:

property (i) means that for every $j,\,1\leqslant j\leqslant q,\,$ there exists a limit direction $s_{\alpha_{\nu(j)}}$ in $P(\Lambda)$ such that its complex conjugate $\overline s_{\alpha_{\nu(j)}}$ is located "near" the direction of the ray $L_{\beta j}$;

wherein one direction $s_{\alpha_\nu}$ may  "be responsible" for several $L_{\beta j}$ (Let usually for all $L_{\beta j}$ and in such a case $\tau=1$); and

property (ii) means that for every $j$ all the straight lines through every $\mu_l\in L_{\beta j}$ with normal direction $\exp(-i\alpha_{\nu(j)})$ (i.e. defined by the normal equations
\begin{equation*}
\Re z\exp i\alpha_{\nu(j)}=\Re\mu_l\exp i\alpha_{\nu(j)}),
\end{equation*}
do not contain any $\mu_k\in\mathcal M$ for $k\not=l$.

In comparison to \cite{MP} (where $q=1$, or $q=2$ and $\beta_1=\beta_2+\pi$) it is found a simplified proof of sufficiency part of the criterion and, due to an idea of S.G. Merzlyakov, a new approach is applied in the proof of the necessity part that allow to eliminate an additional assumption that was prescribed in \cite{MP}.

\medskip
\noindent{\bf Proof of necessity.}
To begin with suppose that considered interpolation is operable and (i) is valid but (ii) fails.

The example given above after the formulation of interpolation problem demonstrates a case when (i) is valid but (ii) fails and interpolation fails too.

For more general demonstration suppose that there exists $j,\,1\leqslant j\leqslant q,\,$ and ${\nu(j)},\, \left|\beta_j+\alpha_\nu(j)\right|<\pi/2$, but there exists $\mu_l\in L_{\beta j}$ and $\mu_k\in \mathcal M$ such that $\Re\mu_l\exp i\alpha_{\nu(j)} =\Re\mu_k\exp i\alpha_{\nu(j)}$. Thus $\Re(\mu_k-\mu_l)\exp i\alpha_{\nu(j)}=0$.

Consider convolution operator $f(z+\mu_l)=f(z+\mu_k)$, generated by $\varphi(z)=\exp\mu_1z-\exp\mu_kz$ with zero set
\begin{equation*}\Lambda=\{\lambda_n=\frac{2\pi n}{|\mu_k-\mu_1|}\exp(i(\pi/2-\beta),\,n\in\bZ\},\, \exp(i\beta)=\frac{\mu_k-\mu_l}{|\mu_k-\mu_l|}, 0\leqslant\beta<2\pi. \end{equation*}
The set $\Lambda$ has two natural limit directions to infinity. Since $\Re(\mu_k-\mu_l)\lambda_n=\Re (2\pi ni)=0$ one of the directions is equal to $ \alpha_{\nu(j)}$ so (i) is valid. We reach a contradiction considering that all the functions $u\in\Sigma(\bC,\Lambda)$ satisfy $f(\mu_l)=f(\mu_k)$ therefore simple interpolation by the functions from $\Sigma(\bC,\Lambda)$ fails for arbitrary interpolation data.

Further suppose that (i) fails and this implies that there exists $j,\,l\leqslant j\leqslant q,\,$ such that none of the directions $s_\alpha,\,$ for $\alpha\in(\beta_j-\pi/2,\beta_j+\pi/2),$ belongs to $P(\Lambda)$.

Using the transformations $z\longmapsto z\exp {(-i\beta_j)}$ and $z\longmapsto z-a_j$ of the complex plane in the setting of the problem of interpolation one readily obtains that without loss of generality it can believed that $\beta_j=0,\,L_j=\bR^+, $ and the set $\bS\cap \{z:\Re z>0\}$ does not contains any direction from $P(\Lambda)$. Then for every $\alpha\in (0,\pi/2)$ the set $\Lambda\cap \{z:\left| \Re z \right|<\left|\Im z\right|\}$ is bounded or empty.

 Consider arbitrary function from $\Sigma(\bC,\Lambda)$:
\begin{equation*} u(z)=\sum_{n=1}^{\infty} c_{n} e^{\lambda_nz},\,z\in \bC,\,\{\lambda_n\}\subset \Lambda.
 \end{equation*}
 The series is absolutely convergent for all $z \in\bC$. As mentioned above the numerical series of norms is convergent therefore applying this for the closed disc $\{z:|z|\leqslant \varepsilon\}$ the numerical series $\sum_{n=1}^{\infty} c_{n} e^{\varepsilon|\lambda_n|}$ is the convergent one.

Suppose the set $\Lambda\cap \{z:\Re z>0\}$ is unbounded.  This implies that $\Re \lambda_{n_k}/ \left| \lambda_{n_k} \right|\to 0, k\to \infty$, for arbitrary subsequence $\{\lambda_{n_k}\}\subset\Lambda,\, \left|\lambda_{n_k}\right|\to\infty$.

 For $x>0$ consider the function $h(x) =\sup\{\Re \lambda_nx-\left|\lambda_n\right|\}$. It is the finite convex function. 
 
The estimate $\left|u(x)\right|\leqslant Ce^{h(\varepsilon x)},\, C>0$, is valid for all $ x>0$.
Indeed,
\begin{equation*}
 \left|u(x)\right|\leqslant \sum_{n=1}^{\infty} \left|c_{n}\right| e^{x\Re\lambda_n}\leqslant e^{h(\varepsilon x)}\sum_{n=1}^{\infty} \left|c_{n}\right| e^{\varepsilon\left|\lambda_n\right|}.
 \end{equation*} 
Simple interpolation $f(\mu_k)=b_k,\, \mu_k >0$ fails, e.g. for $b_k\geqslant k\exp h(\varepsilon\mu_k$.

In the case when the set $\Lambda\cap \{z:\Re z>0\}$ is bounded or empty it may be noted that $h(\varepsilon x)\leqslant Cx$ or $h(\varepsilon x)\leqslant 0$, respectively, and interpolation fails too. Proof of necessity is complete.

Thus in the proof it is disclosed a method of estimation of growth on the rays of modulus of sum of the series of exponentials that have the specific location of limit directions at infinity for $\Lambda$. That may allows also to treat multiple interpolation problems mentioned above.

Let us consider the countable set of real parts of all $\lambda_n\in\{z:\Re z>0\}$ then put it in order of increase and denote it as $\{l_k\}$. Determine the complex sequence $\widetilde\lambda_k=l_k+i\min\{\left| \Im \lambda_n\right|: \Re \lambda_n=l_k>0\}$. Further let us define the convex (in extended sense) function $\Phi:\bR^+\rightarrow [0,+\infty]$ as follows:
\begin{equation*}
\Phi(t)=\sup\{at+b:al_k+b\leqslant\left | \widetilde\lambda_k\right |,\,k\in \bN\}. 
\end{equation*}
Without difficulty one obtains that $h(x)=\Phi^*(x)=\sup_{t>0}\{xt-\Phi(t)\}$. That is the Young-Fenchel transform is operable to estimate $u(x)$. 

For a time the proof of sufficiency is postponed and, to begin with, the problem is reduced to the interpolation problem in the kernel of specific convolution operator  and then observations of duality using the Laplace transformation of functionals are applied. As the result, the proof of existence of representations (\ref{F3}) may be reduced to the consideration of equivalent natural dual statements in the classical space $P_{\bC}$ of all entire functions of exponential type.

\noindent{ \bf Reduction to the problem of interpolation by functions from the kernel of a specific convolution operator.}
It should be noted the following. As soon as a subset $\Lambda_1\subset\Lambda$ is found so as to establish the representation $H(\bC)=\Sigma (\bC,\Lambda_1)+I(\bC,\mathcal M)$  for some $\mathcal M\subset \bC$, the representation (\ref{F1}) is also valid, i. e. it is solvable the considered problem of interpolation by functions from $\Sigma (\bC,\Lambda)$ with the node set $\mathcal M$.

Let us consider the specific procedure of extraction of a sparse subsequence from $\Lambda$  while retaining the given set of its limit directions at infinity. 
\begin{lemma}\label{L1}
Let $\Lambda$ be an arbitrary unbounded set and $\{s_{\alpha_1}, s_{\alpha_2},\cdots, s_{\alpha_{\tau}}\}\subset \Lambda$. There exists a sequence $\{\lambda_n\}\subset \Lambda$ that possesses the properties:

1. $P(\{\lambda_n\})=\{s_{\alpha_1}, s_{\alpha_2},\cdots,s_{\alpha_{\tau}}\}$.

2. $|\lambda_{n+1}|>2|\lambda_{n}|$.
\end{lemma}
Proof.  Take arbitrary $\lambda_1\in \Lambda,\,$ and then select, using definition of limit directions at infinity,  by induction $\lambda_2,\, |\lambda_2/|\lambda_2|-s_{\alpha_2}|<1/2^2,\,|\lambda_2|>2|\lambda_2|$, ...  ,$\lambda_{\tau},\, |\lambda_{\tau}/|\lambda_{\tau}|-s_{\alpha_{\tau}}|<1/2^{\tau},\,|\lambda_{\tau}|>2|\lambda_{\tau-1}|,\,\lambda_{\tau+1},\, |\lambda_{\tau+1}/|\lambda_{\tau+1}|-s_{\alpha_{\tau+1}}|<1/2^{\tau+1},\,|\lambda_{\tau+1}|>2|\lambda_{\tau}|,\cdots$  This selection completes the proof.
  
Denote $\Lambda_1=\{\lambda_n\}$. For arbitrary set $\Lambda$ the subspace $\Sigma(\Lambda, \bC)$ used in (\ref{F1}) is generally not closed in $H(\bC)$.   It will be demonstrated thereafter that $\Sigma(\Lambda_1, \bC)$ coincides with the kernel of some continuous convolution operator defined in $H(\bC)$, so $\Sigma(\Lambda_1, \bC)$ is the closed subspace of $H(\bC)$.  

The sequence $\Lambda_1$ has density zero. Consider  an entire function which have simple zeroes $\lambda_{n}$, 
 \begin{equation*}
 G(z)=\prod_{n=1}^{\infty} \left( 1-\dfrac{z}{\lambda_{n}} \right).
 \end{equation*}
 The function $ G $ is exponential type function and it has minimal type of growth.

 The quantity
$$\delta= \limsup_{n\to\infty}\dfrac{1}{|\lambda_{n}|}\ln\dfrac{1}{\bigl|G^\prime(\lambda_{n})\bigl|}$$ 
 is the condensation index of Gelfond-Leont'ev. The sparse structure of $\Lambda=\{\lambda_n\}$ implies (\cite{MP}) that $ \delta= 0$.
 
From the results of monograph of A.F. Leont'ev \cite{Leon} (Theorem 4.2.4)  it follows the important proposition. 

\medskip
{\noindent \it  Let $ \delta=0$ and the set $\{\lambda_n\}$ has the density zero. Then any function in the closure, in the topology of $ H(\bC)$,  of linear hull of the system $\{\exp\lambda_nz\}$  is represented as a series of exponentials.}

The function $G$ generates linear continuous surjective convolution operator $M_{G}$ on the space $H(\bC)$. The kernel of the operator $\Ker M_{G}$ is a closed subspace which is invariant in respect to operator of differentiation and admit the spectral synthesis \cite{Leon} (Theorem 4.2.9). i.e. it coincides with the closure of the linear hull of all exponentials $e^{\lambda_nz}\in \Ker M_{G}$. 
Then, taking into account \cite{Leon} (Theorem 2.1.3), one obtains the following proposition. 
\begin{lemma} \label{L2} The kernel $\Ker M_{G}$ coincides with the set of all entire functions $f(z)$, that are represented by the series of exponentials ,
\begin{equation*}
f(z)=\sum_{n=1}^\infty c_{n} e^{\lambda_{n}z},\,z\in \bC,\end{equation*}
converging in the topology of the space $H(\bC)$, that is $ \Ker M_{G}= \Sigma(\Lambda_1, \bC)$.
\end{lemma}
The assertion of lemma~\ref{L2} is referred to as the Fundamental Principle. It should be noted that the Principle is studied for one and several complex variables and in \cite{Kri1}, \cite{Kri2}, \cite{Kri4} a new local condensation quantity $S$ is introduced and criteria are obtained. The results give an alternative proof of the lemma.

 In view of lemma~\ref{L2}, for a given arbitrary unbounded set $\Lambda$, solubility of the problem of interpolation by functions from $\Sigma(\bC, \Lambda)$  with some node set $\mathcal M$ is a consequence of representation
\begin{equation}\label{F3}
H(\bC)=\Ker M_G + I(\bC,\mathcal M).
\end{equation}
In fact these representations are matched with the sequence $ \{\lambda_n\}=\Lambda_1\subset \Lambda$, so they imply representations (\ref{F1}).  Representations (\ref{F3}) will be considered in the proof of sufficiency of criterion of interpolation. 

\medskip
\noindent{ \bf Dual formulation of the interpolation problem.}
In the following it is used the methods of proof (see \cite{NapP} and \cite{MP}, \cite{MP1}) that are different from that of V.V. Napalkov . The proof is based on the duality using the Laplace transformation $\mathcal L$ of functionals and representations (\ref{F3}) are equivalent to natural dual statements in the space $P_{\bC}$.

Denote by $P_{\bC}$ the space of all entire functions of exponential type. It is introduced in $P_{\bC}$ a traditional inductive topology of weighted Banach spaces of entire functions in virtue of which topological isomorphism is provided  between the strong dual $H^{*}(\bC)$ and the space  $P_{\bC}$, implemented by the Laplace transformation $\mathcal L$ of functionals $F\in H^{*}(\bC)$. To be more precise, the linear continuous one-to-one Laplace transformation $\mathcal L$ of functionals $F\in H^{*}(\bC)$ is defined as follows: $ \mathcal L: F\longmapsto \mathcal L{F}(z)=\left\langle F_\lambda,e^{\lambda z}\right\rangle,\, \mathcal L{F}\in P_{\bC}$.
  
The inductive topology in ${(LN}^{*})$-space $P_{\bC}$ is not described in the terms of sequence convergence, however all sequentially closed subspaces are the closed ones (\cite{Seb}). The exact definition of sequence convergence will be formulated later in the proof of sufficiency of the main theorem.

Separately continuous bilinear form $\left[\cdot, \cdot\right]: H(\bC)\times P_{\bC}\longmapsto \mathbb C$, defined according to the formula $\left[ \psi, \varphi\right]=\left<{\mathcal L}^{-1}\varphi,\,\psi\right>,\,\psi\in H(\bC),\,\varphi \in P_{\bC}$, gives the duality $\left[H(\bC), P_{\bC}\right]$.
 The mapping $\varphi \longmapsto \left[\cdot,\phi \right]= \left<{\mathcal L}^{-1}\varphi,\cdot \right>$, where $ {\mathcal L}^{-1}\varphi  \in H^{*}(\bC)$, establishes an isomorphism between $P_{\bC}$ and strong dual $H^{*}(\bC)$. According to the defined duality, every function in the space ${P_{\bC}}$ is in one-to-one correspondence with the unique linear continuous functional from $H^*(\bC)$.
 
Function $ G \in P_{\bC},\,G\not\equiv0$, having minimal type of growth, generates in the space $H(\bC)$ convolution operator $M_G: H(\bC)\longmapsto H(\bC)$, which may be defined in considered duality as
\begin{equation*} M_G [\psi](z)=\left[S_z\bigl(\psi(\lambda)\bigr),  G_\lambda\right]=\left<({\mathcal L}^{-1}G)_\lambda, \psi (z+\lambda)\right>,
\end{equation*}
wherein $ S_z $ is the shift operator: $  S_z\bigl(\psi(\lambda)\bigr)=\psi(\lambda+z)$. 

Denote by $\Ker M_{G}=\{f\in H(\bC): M_G[f]=0\}$ the kernel of convolution operator $ M_{G}$. It is the closed subspace of $H(\bC)$ and is invariant in respect to the operator of differentiation and it admits the spectral synthesis.  

For any $\mathcal M$ there exists a function $\psi=\psi_{\mathcal M}\in H(\bC)$ such that its zero set coincides with all $\mu_k\in \mathcal M$ with multiplicities $m_k$. It may be readily observed that ideal $I(\bC,\mathcal M)$ is equal to the closed principal ideal $\bigl(\psi\bigr)$ generated by the function $\psi$:
\begin{equation*}
I(\bC,\mathcal M)=\bigl(\psi\bigr)=\{h\in H({\bC}): h=\psi\cdot r, \, r\in H({\bC})\}.
\end{equation*}
For any node set $\mathcal M$ the existence of representation (\ref{F3}) is equivalent to the following two statements.

\medskip
 \noindent{\it  $\it (I)$ Subspace $\Ker M_{G}+ \bigl(\psi\bigr)$ is the dense subspace in  $H(\bC);$}  

\noindent {\it $\it (II)$ Subspace $\Ker M_{G}+ \bigl(\psi\bigr)$ is the closed subspace of $H(\bC)$. }

If $U$ is the subspace of topological vector space $X$ denote by $ U^\circ$ its polar set (or annihilator) that is the set of all functionals in $X^*$ which is zero on $U$. 

The statement {\it  $\it (I)$} is equivalent to the statement that $\bigl(\Ker M_{G}+ (\psi)\bigr)^\circ=\bigl(\Ker M_{G}\bigr)^\circ\cap\bigl((\psi)\bigr)^\circ=\{0\}$.
$H(D)$ is the (FS)-space. 
In the case (FS) (or (DFS)) spaces Lemma 2 from the paper \cite{Mer1} (see also \cite{NapP} for detail) implies that the statement {\it  $\it (II)$}  is equivalent to the statement that the subspace  $\bigl(\Ker M_{G}\bigr)^\circ+ \bigl((\psi)\bigr)^\circ$  is closed in $P_{\bC}$ which is the (DFS)-space.  

Thereafter it is discussed the duality implementation of the annihilators in the space $P_{\bC}$, that is the ring.

Adjoint operator to the convolution operator $M_{G}$ is the operator $ A_G $ of multiplication on characteristic function $G$ which is correctly defined on the functions ${v \in P_{\bC}}$ in the following manner: $ v \longmapsto G\cdot v$. 
 According to the duality, annihilator $\bigl(\Ker M_{G}\bigr)^\circ$ coincides with the closed principal ideal which is defined as
\begin{equation*}
\bigl(G\bigr)_{P_{\bC}}=\{h\in P_{\bC}:  h= G\cdot v;\,v\in P_{\bC}\}.
\end{equation*}
This is the range of the adjoint operator $A_G$ in $P_{\bC}$. 

The Lindel\"of theorem implies the theorem of division in $P_{\bC}$ on the function $G$ that has minimal type of growth and it gives (see \cite{MP}) that $\bigl(G\bigr)_{P_{\bC}}=\bigl(G\bigr)\cap {P_{\bC}}$ where $\bigl(G\bigr)$ is the principal ideal in $H(\bC)$ generated by the function $G$. In particular this implies that the ideal $\bigl(G\bigr)_{P_{\bC}}$ is closed  because the topology of the space $P_{\bC}$ is stronger than the topology of point-wise convergence.

It is known that $(M^{*})$-space $H(\bC)$ is reflexive so its strong second dual $H^{**}(\bC) $ is canonically isomorphic to the $H(\bC)$. Therefore the mapping $\psi \longmapsto \left[\psi, \cdot\right]$, in view of this  canonical isomorphism, defines the isomorphism between $(M^{*})$-space $H(\bC)$ and the strong dual $P_{\bC}^{*}$. Every function in $H(\bC)$ is in one-to-one correspondence with the unique linear continuous functional in the strong dual $P_{\bC}^*$. 

Every function $\psi\in H(\bC),\,\psi\not\equiv0, $ generates in the space $P_{\bC}$ convolution operator $\widetilde M_\psi: P_{\bC}\longmapsto P_{\bC}$, (see \cite{MP} for more detail)

\begin{equation*} \widetilde M_\psi[\varphi](z)=\left[(\Theta\psi)_\lambda, S_z\bigl(\varphi(\lambda)\bigr)\right],\,\varphi\in P_{\bC}=\bigl<{ (\mathcal L}^{-1}\varphi)_\lambda,\,e^{z\lambda}\psi(\lambda)\bigr>,
\end{equation*}
wherein $ S_z $ is the shift operator, $ S_z\bigl(\varphi(\lambda)\bigr)=\varphi(\lambda+z),\,\lambda \in\mathbb C$.

 Furthermore using the integral representation of the Borel transform (\cite{Leon}) one readily obtain that the convolution operator $\widetilde M_\psi$ has the following representation:
 \begin{equation*}
 \widetilde {M}_{\psi}[\varphi](z)=\frac{1}{2\pi i}\int \limits_{C}\psi(\lambda)e^{z\lambda}\gamma_\varphi (\lambda)\,d\lambda,\,\varphi\in P_{\bC},
 \end{equation*}
  where $\gamma_\varphi$ is the Borel transform of function $\varphi$, and $C$ is rectifiable closed curve which surrounds all singularities of the function $\gamma_\varphi$.
It is known \cite{10} that $\widetilde M_\psi$  is the linear continuous surjective operator that 
is adjoint operator for the operator $\widetilde A_{\psi}$ of multiplication on the function $\psi$ in the space $H(\bC)$. The operator $\widetilde A_{\psi}$ is defined on the functions $g\in H(\bC)$ in the following manner: $ g \longmapsto \psi\cdot g$ and it is linear and continuous and his range coincides with the closed ideal $(\psi)$.  Denote $\Ker  \widetilde M_{\psi}=\{f\in P_{\bC}: \widetilde M_\psi[f]=0\}$. According to the duality $\bigl( (\psi)\bigr)^\circ$ is equal to $\Ker \widetilde M_{\psi}$. This is closed subspace of $P_{\bC}$ that admits the spectral synthesis. 

 Moreover  the Fundamental Principle for $\Ker \widetilde M_{\psi}$ in $P_{\bC}$  is valid (see e.g.  \cite{10}). To be more precise this means the following needful statement:
 
\medskip
\noindent {\it Subspace $\Ker \widetilde M_{\psi}\subset P_{\bC}$ is equal to the linear hull of all monomes $\{z^\nu e^{\mu_{k} z}\},\, k\in\mathbb N,\,\nu=0, 1, \cdots,m_{k}-1$ thus it consist of all polynomials of exponentials of the following form} 
\begin{equation}\label{F4}
p(z)=\sum_{\Fin_{p}(\mathcal M)} a_k(z)e^{\mu_k z},
\end{equation}
{\noindent\it wherein $ a_k(z) $ is some polynomials of degrees not exceeding $m_k-1$. On the right side there is the sum over all $\mu_k$ in some finite set $\Fin_{p}(\mathcal M)\subset \mathcal M$.} 

Thus it is obtained on the base of considered duality that statements {\it $\it (I)$} и {\it $\it (II)$} in $(M^{*})$-space $H(\bC)$ are equivalent to the following two corresponding  statements in $({LN}^{*})$-space ${P_{\bC}}$ : 

\medskip
\noindent {\it  $\it (I^{*}) $ The equality ${\bigl(G\bigr)_{P_{\bC}}}\cap \Ker \widetilde M_{\psi}=\{0\}$ is valid.}

\noindent{\it $(\it II^{*})$  The subspace ${\bigl(G\bigr)_{P_{\bC}}}+\Ker \widetilde M_{\psi}$ is the closed subspace of ${P_{\bC}}$.}

As a result we obtain the following.  
\begin{theorem}\label{T2} Existence of representations (\ref{F3}) of the space $H(\bC)$ is equivalent to the dual statements $\it (I^{*})$ и $(\it II^{*})$ in the space ${P_{\bC}}$.
\end{theorem}
The theorem may be proved under less restrictive assumptions as stated in \cite{NapP}.  
\section{Proof of sufficiency}
Lemma~\ref{L1} allows to extract an sparse sequence $\Lambda_1=\{\lambda_n\}$ from $\Lambda$.  If it is not explicitly mentioned thereafter it will be supposed in the proof of sufficiency that $\Lambda$ is equal to the sequence $\{\lambda_n\}$ and $ \Ker M_{G}= \Sigma(\Lambda, \bC)$. Thus according to foregoing dual formulation of the interpolation problem natural dual statements $\it (I^{*})$ and $(\it II^{*})$ should be proved in the space  ${P_{\bC}}$.

Let us prove the statement $(\it II^{*})$. It is well-known (\cite{Seb}) that in $({LN}^{*})$-space $P_{\bC}$ any subspace $X$ is closed  if and only if it is sequentially closed.
 
 Consider arbitrary sequence $\{g_l\}_{l\in\mathbb N}$ of functions in subspace $(G)_{P_{\bC}}+\Ker \widetilde M_{\psi}$. Suppose it converges in the space $P_{\bC}$ to the function $g \in {P_{\bC}}$. We need to show that the limit function $g$ belongs to $(G)_{P_{\bC}}+\Ker \widetilde M_{\psi}$.
 
It should be noted, in $(LN^{*})$-topology of the space $P_{\bC}$ convergence of a sequence $\{g_l\}_{l\in \mathbb N}$ to the limit function $ g $  is equivalent to the following two assertions: 

\noindent 1. The sequence $\{g_l\}$ converges to $ g $ in the topology of the space  $ H(\mathbb C)$, thus it  converges uniformly on all compact sets in $\bC$.

\noindent 2. There exists a constant $A>0$ and a convex compact set $K\subset \bC$, such that for all $ l\in\mathbb N $ the following estimate is satisfied:
\begin{equation}\label{F5}
|g_l(z)|\leqslant Ae^{H_K(z)},\, z\in\mathbb C.
\end{equation}
 Here $ H_K(z)=\sup_{\sigma\in K} \Re (z\sigma)$ is the support function of compact set $K$ in the sense of complex analysis. The function $H_K(\theta)$ is continuous and positive homogeneous. If $ z=|z|e^{i\theta},\ $  $H_K(z)=h_K(\theta)|z|$, where $h_K(\theta)=\sup_{\sigma\in K} \Re (\sigma e^{i\theta})$ is the support function (in the sense of $\bR^2$) of the compact set $\overline K$ which is complex conjugate of $K$. The function $h_K(\theta)$ is continuous.
 
Sometimes it is appropriate to interpret the value $h_K(\theta)$ as supremum of projections on the ray $\{te^{-i\theta},\,t>0\}$ over all points of compact set $K$. By other words, denote by $ k(\theta, K): \mathbb C\mapsto \bR$ the support function (in the sense of $\bR^2$) of the compact set $K$ in the direction $e^{i\theta} =(\cos\theta, \sin\theta)$,  then $k(\theta, K)=h_K(-\theta)$.  It may be noted, that the support half-plane $\Pi(\overline{s}, K)$ (in the sense of $\bR^2$) of convex compact set $K$ with the normal direction $\overline{s}_\theta=e^{-i_\theta}$ has the form $\Pi(\overline{s}, K)=\Pi_0(\overline{s}_\theta)+\overline {s}_\theta h_K(\theta)$, where $h_K(\theta)=k(-\theta, K)$. Here $ \Pi_0(\overline{s}_\theta)=\{z=x+iy\in\mathbb C: \Re(s_\theta z)=x\cos(-\theta)+y\sin(-\theta)<0\})$. Convex compact set $K$ is equal to the intersection of all closed half-planes $\overline\Pi(\overline{s}, K)$ over all $s\in \bS$.

The sequence  $\{g_l\}$ has the form $ g_l= p_l+R_l, $ where for every $l\in\bN$ the function $p_l\in\Ker \widetilde M_{\psi}$ and the function $ R_l \in  (G)_{P_{\bC}}$, in particular $ R_l(\lambda_n)=0,\, n\in\bN$.
  
In the case when the sequence $\{g_l\}$ contains infinite number of elements with $ R_l \equiv 0$ the limit function $ g$ belongs to $\Ker \widetilde M_{\psi}$, due to continuity of the convolution operator.  When the sequence $\{g_l\}$ contains infinite number of elements with $ p_l\equiv 0$ the limit function $ g$ belongs to $(G)_{P_{\bC}}$, as the topology in $P_{\bC}$ is stronger than  topology of point-wise convergence and due to the theorem of division in $P_{\bC}$ on the function $G$ that has minimal type of growth.

It is obtained for all sequences $\{g_l\}$ from both categories that $g\in(G)_{P_{\bC}}+\Ker \widetilde M_{\psi}$. Therefore without loss of generality it may believed thereafter that $\{g_l\}=\{p_l+R_l\}$ is such that $R_l\not \equiv 0,\,p_l\not\equiv 0$, for all $ l\in\bN$.

For arbitrary fixed $l\in\bN$ since $p_l \in \Ker  \widetilde M_{\psi},\,  p_l\not\equiv 0, $ then due to the Fundamental Principle (\ref{F4}) it is the polynomial of exponentials of the form
\begin{equation}\label{F6}
p_l(z)=\sum_{\Fin_{p_l}(\mathcal M)} a^l_k(z)e^{\mu_k z},
\end{equation}
wherein $ a_k(z) $ is some polynomials of degrees not exceeding $m_k-1$. On the right-hand side of the determination the sum stands which is over all $\mu_k$ in some finite set $\Fin_{p}(\mathcal M)\subset \mathcal M$. 

Let us use the disc K=$\{z:|z|\leqslant B\}$ as the compact set in description (\ref{F5}) of convergence in ${P_{\bC}}$. Convergence  of $\{g_l\}$ in $P_{\bC}$ implies that 
$|p_l(z)+R_l(z)|\leqslant Ae^{H_K(z)},\, z\in\mathbb C$, but $R_l(\lambda_n)=0$, so that
\begin{equation}\label{F7}
|g_l(\lambda_n)|=|p_l(\lambda_n)|\leqslant Ae^{B\left |\lambda_n\right|)}
\end{equation}
for all $n\in \bN$.

Increasing if required the value of $B$ one may assume, without any loss of generality, that the closed disc $K$ contains all points $\mu_k$  lying in finite number outside $\cap_{j=1}^q F_j,\,$ where $F_j=F\cup L_{\beta_j},$ and all initial points $a_j$ of the rays $L_{\beta_j}=\{a_j+t\exp(i\beta_j),\, t>0\},\,j=1,\cdots,q$. 

It will be established thereafter that for every $l\in\bN$ all exponents $\mu_k\in\Fin_{p_l}(\mathcal M)$ in representation (\ref{F6}) belong to some compact set $K_1$.

If $\Fin_{p_l}(\mathcal M)\subset K$ all is done. Suppose there exists a point $\mu_k\in\Fin_{p_l}(\mathcal M)$ outside $K$ then $\mu_k\in L_{\beta_j}$ for some $j$. By property (i) there exists $s_{\alpha_{\nu(j)}}\in P(\Lambda)$, such that   $|\beta_j+\alpha_\nu(j)|<\pi/2$. 

To begin with consider all $\mu_k$ that additionally satisfy $Re\mu_k\exp {i\alpha_{\nu(j)}}>B$ and denote by $\mu_\tau$ the node for which $\max Re\mu_k\exp {i\alpha_{\nu(j)}}$ over all $\mu_k\in\Fin_{p_l}(\mathcal M)$ is attained. Then by property (ii) the unique point  $\mu_\tau$ is situated on the support line $\Re z\exp i\alpha_{\nu(j)}=\Re\mu_\tau\exp i\alpha_{\nu(j)}$.  

Denote by $T$ the closed convex hull of the set $\Fin_{p_l}(\mathcal M)\cup K$. Compact set $T$ coincides with the closed convex hull of its extreme points. By the foregoing reasoning $T\not=K$ therefore the point $\mu_\tau\in\Fin_{p_l}(\mathcal M),\,$ $\mu_\tau\not\in K,$ is the extreme point of compact set $T$. 

 Compact set $T$ is equal to the intersection of all closed half-planes $\overline\Pi(\overline{s}, K)$ over all $s\in \bS$.
 
%Hence there exists such a direction $s_{\alpha_{\nu(j)}}=e^{i\alpha_{\nu(j)}}$, that on the boundary of the support half-plane $\Pi(\overline s_{\alpha_{\nu(j)}},T)$  is situated only one extreme point $\mu_\tau\in\Fin_{p_l}(\mathcal M)$.
 The boundary of $\Pi(\overline s_{\alpha_{\nu(j)}}, T)$ is the support line $\{z=x+iy:\Re ze^{i\alpha_{\nu(j)}}=x\cos(-\alpha_{\nu(j)})+y\sin(-\alpha_{\nu(j)})=h_T(\alpha_{\nu(j)})\}$ with the normal direction
 $\overline s_{\alpha_{\nu(j)}}=e^{-i\alpha_{\nu(j)}}=(\cos(-\alpha_{\nu(j)}), \sin(-\alpha_{\nu(j)}))$. Note also that $h_T(\alpha_{\nu(j)})=\Re\mu_\tau e^{i\alpha_{\nu(j)}}$.

Denote $r_l(z)=p_l(z)-a^q_k(z)e^{\mu_\tau z},\,$ $\Fin_{r_l}(\mathcal M)=\{\mu_k\in\Fin_{p_l}(\mathcal M):\mu_k\not=\mu_\tau\},\,$ and denote by $T_1$ the closed convex hull of the set $\Fin_{r_l}(\mathcal M)\cup K$. Then by assumption $T\not=T_1$. Note that $K\subset T_1$ and occasionally $K=T_1$.

Actually in virtue of the definition of the direction $s_{\alpha_{\nu(j)}}$ in (i) and by the property (ii), compact set $\Fin_{r_l}(\mathcal M)\cup K$ is situated in the open support half-plane $\Pi(\overline s_{\alpha_{\nu(j)}}, T)$. That means that $\Re z e^{i\alpha_{\nu(j)}}<\Re\mu_\tau e^{i\alpha_{\nu(j)}}$ for all $z\in\Fin_{r_l}(\mathcal M)\cup   K$. By continuity of $\Re z e^{i\alpha_{\nu(j)}}$ there exists $z_0\in \Fin_{r_l}(\mathcal M)\cup K$ such that
\begin{equation*}
\sup_{z\in\Fin_{r_l}(\mathcal M)}\Re z e^{i\alpha_{\nu(j)}}=\Re z_0 e^{i\alpha_{\nu(j)}}<\Re\mu_\tau e^{i\alpha_{\nu(j)}}.
\end{equation*}
Further, since $\Re z_0 e^{i\alpha_{\nu(j)}}= h_{T_1}(\alpha_{\nu(j)})$ we get that for some $\varepsilon>0$ the following inequality is valid:
 $\Re\mu_\tau e^{i\alpha_{\nu(j)}}>h_{T_1}(\alpha_{\nu(j)})+2\varepsilon$. By continuity it follows  that for some $\delta >0,\,$ $\Re\mu_\tau e^{i\theta}>h_{T_1}(\theta)+2\varepsilon,\, |\theta-\alpha_{\nu(j)}|<\delta$.

Further it is used the estimate and the fact that for all $z=re^{i\theta}$ polynomial of exponentials  $r_l(z)=p_l(z)-a^q_k(z)e^{\mu_\tau z}$ satisfies the estimate $|r_l(z)|\leqslant C_\varepsilon e^{(h_{T_1}(\theta)+\varepsilon)|z|}$. The estimate is a consequence of the integral representation of the Borel transform.

For all $z=re^{i\theta}$ in the angle $\Gamma_\delta(\alpha_{\nu(j)})=\{z=re^{i\theta}:  |\theta+\alpha_{\nu(j)}|<\delta\}$ for $\{|z|>r\}$ (in particular, the disc $\{|z|\leqslant r\}$ may contain all eventual common zeroes of polynomials $ a^l_k(z) $ in (\ref{F6})), from the estimate $ |p_l(z)|\geqslant |a^q_k(z)e^{\mu_\tau z}|-|r_l(z)|$ it follows that 
 \begin{equation*}
 |p_l(z)|\geqslant e^{\Re\mu_\tau e^{i\theta}|z|}( 1-C_\varepsilon e^{(h_{T_1}(\theta)+\varepsilon)|z|-\Re\mu_\tau e^{i\theta}|z|})> e^{\Re\mu_\tau e^{i\theta}|z|}( 1-C_\varepsilon e^{-\varepsilon |z|}).
\end{equation*}
The values of $r$ and some of the constants are depending on $l$ but this is irrelevant in the setting.

In virtue of the choice of the direction $s_{\alpha_{\nu(j)}},\,$ $\Re\mu_\tau e^{i\alpha_{\nu(j)}}=h_T(\alpha_{\nu(j)})$.

One may find $\delta_1 >0,\,\delta_1<\delta$, in such a way that $\Re\mu_\tau e^{i\theta}=h_T(\theta)$ for all $|\theta-\alpha_{\nu(j)}|<\delta_1$. Indeed, the point $\mu_\tau$ is the common point of two adjacent segments on the boundary of $T$ and it is sufficient to take any $\delta$ that is less than the value of acute angle between normal directions to the adjacent segments.

It is obtained that for all $z$ in the angle $\Gamma_{\delta_1}(\alpha_{\nu(j)})\subset\Gamma_{\delta}(\alpha_{\nu(j)}),\,$ the following estimate is valid, $\Re\mu_\tau e^{i\theta}|z|=H_T(\theta)>H_{T_1}(\theta)+2\varepsilon|z|$.

Furthermore the foregoing estimate of $ |p_l(z)|$ implies the estimate $ |p_l(z)|>Ae^{H_{T_1}(z)+\varepsilon|z|}$, for all $z=re^{i\theta}$ in some set $\Gamma_{\delta_1,r_1}(\alpha_{\nu(j)})= \Gamma_{\delta_1}(\alpha_{\nu(j)})\cap \{|z|>r_1\},\,r_1>r$, wherein $A$ is the constant in inequality (\ref{F7}). 

An infinite sequence of exponents $\{\lambda_{n_k}\}\subset\Lambda$ is situated in the set $\Gamma_{\delta_1,r_1}(\alpha_{\nu(j)})$ because the direction  $s_0=e^{i\alpha_{\nu(j)}}$ is a limit direction at infinity for $\Lambda$ by the assumptions of the theorem~\ref{T1}.

Thus the following estimate is obtained,
\begin{equation*}
|p_l(\lambda_{n_k})|>Ae^{H_{T_1}(\lambda_{n_k})+\varepsilon|\lambda_{n_k}|},\, \{\lambda_{n_k}\}\subset\Gamma_{\delta_1,r_1}(\alpha_{\nu(j)}),
\end{equation*} 
that is contrary to the estimate (7) because $H_K(z)\leqslant H_{T_1}(z)$ for all $z$.

We obtain that $\Fin_{p_l}(\mathcal M)$ does not contains considered point $\mu_\tau$ outside $K$. Repeating if required the reasoning above we obtain that the set  $\Fin_{p_l}(\mathcal M)$ does not contains any point lying outside $K$ and satisfying  $\Re\mu_k\exp {i\alpha_{\nu(j)}}>B$.   

Further suppose there exists a point $\mu_k\in\Fin_{p_l}(\mathcal M)$ outside $K$ satisfying $\Re\mu_k\exp {i\alpha_{\nu(j)}}\leqslant B$ wherein $\mu_k\in L_{\beta_j}$ for some $j$ and, by property (i), there exists the matching direction $s_{\alpha_{\nu(j)}}\in P(\Lambda)$, such that   $|\beta_j+\alpha_\nu(j)|<\pi/2$.

Consider all $\mu_k$ as described. Recall that all $a_j\in K$. Fix $j$ and note that all the points $\mu_k$ on $L_{\beta_j}$ have the following form $\mu_k=a_j+t_k\exp {i\beta_j}$. Hence
\begin{equation*}
\Re a_j\exp {i\alpha_{\nu(j)}}+\Re t_k\exp i(\beta_j+\alpha_{\nu(j)})\leqslant B
\end{equation*}
and further $t_k\cos(\beta_j+\alpha_{\nu(j)})\leqslant B-\Re a_j\exp {i\alpha_{\nu(j)}}$,
where $\cos(\beta_j+\alpha_{\nu(j)})>0$ by the property (ii) that provides that none $L_j$ is parallel to the support line of $K$ with normal direction $\overline s_{\alpha_{\nu(j)}}$. Denote 
\begin{equation*}
T_j=\frac{B-\Re a_j\exp {i\alpha_{\nu(j)}}}{\cos(\beta_j+\alpha_{\nu(j)})}.
\end{equation*}

 It is shown that for every $j$ all points $\mu_k\in\Fin_{p_l}(\mathcal M)$ in the considered category are  situated on the segments of $L_j$ of the form $a_j+t\exp i\beta_j,\,0\leqslant t\leqslant T_j$, wherein $T_j$ is not depending of $l$. 

As a result of the demonstration above we obtain the following: For all $l\in \bN$ sets $\Fin_{p_l}(\mathcal M)$ belongs to the compact set $K_1$, where $K_1$ is the union of $K$ and finite number of the determined segments that are independent of $l$.   

The set $\mathcal M$ is discrete in $\bC$, so it follows that in any sequence of the form $ g_l= p_l+R_l,\,R_l\not \equiv 0,\,p_l\not\equiv 0, $ that is convergent in $P_{\bC}$, the sequence $\{p_l\}$ belongs to a finite-dimensional subspace $X\subset \Ker  \widetilde M_{\psi}$.

In topological vector space  the algebraical sum of finite-dimensional and closed subspaces is the closed subspace \cite{Rud}. Therefore the limit function $g$  of the sequence $ g_l= p_l+R_l$  belongs to $\Ker  \widetilde M_{\psi}+(G)_{P_{\bC}}$. The proof of  $(\it II^{*})$ is complete.

The foregoing proof implies statement $(\it I^{*})$. Indeed, to establish $(\it I^{*})$ we need to show that none of polynomials of exponentials $p\in\Ker \widetilde M_{\psi},\,p\not\equiv 0$, may belong to  $\bigl(G\bigr)_{P_{\bC}}$. 
 
Assuming the contrary, let $p\in\bigl(G\bigr)_{P_{\bC}},\,p\not\equiv 0$ then $p(\lambda_n)=0,\,n\in\bN$.
  
Consider the stationary sequence $p_l=p$ for all $l\in\bN$. Polynomial $p$ has the form (\ref{F4}). Select a compact set $\widetilde K\subset \bC$ in such a manner that $\widetilde K$ does not contain any point $\mu_k\in \Fin_{p}(\mathcal M)$. We have an obvious estimate $0=|p_l(\lambda_n)|<Ae^{H_{\widetilde K}(\lambda_n)},\,n\in\bN$. But as shown above in the foregoing proof the estimate implies that all $\mu_k\in \Fin_{p}(\mathcal M)$ should be situated in $\widetilde K_1.$ Thus we reach a contradiction so $p\equiv0$. Statement $(\it I^{*})$ is proved.
 
Certainly the statement  $(\it I^{*})$ may be proved by the direct-acting argument. Let $p$ be an polynomial of exponentials of the form~(\ref{F4}). Let us consider the convex hull $T$ of extreme points of  $\Fin_{p}(\mathcal M)$ (it is a convex polygon) then by reasoning as at the beginning of the foregoing demonstration one readily obtains the simple estimate
\begin{equation*}
|p(z)|>e^{\Re\mu_\tau e^{i\alpha_{\nu(j)}}|z|}( 1-C_\varepsilon e^{-\varepsilon |z|})>0,
\end{equation*}
for all $z,|z|>r$, in an angle $\Gamma_{\delta}(\alpha_{\nu(j)})$, for some direction $s_{\alpha_{\nu(j)}}=e^{i\alpha_{\nu(j)}}$ relevant to arbitrary extreme point $\mu_\tau$ of compact set  $T$. Thus the polynomial $p$ of exponentials does not have any zeroes in some angle $\Gamma_{\delta}(\alpha_{\nu(j)})$ outside a disc $\{z:|z|\leqslant r\}$. But the set should contain infinite number of point from $\Lambda$ because of $s_{\alpha_{\nu(j)}}\in P(\Lambda)$, in virtue of the hypotheses of the theorem. It is obtained that $p\equiv0$, if $p\in\bigl(G\bigr)_{P_{\bC}}$.

Theorem~\ref{T2} states existence of representation (\ref{F3}) for $\Lambda_1=\{\lambda_n\}$ from Lemma~\ref{L1}, wherein $\{\lambda_n\}\subset \Lambda$, and this gives us that the representation (\ref{F1}) for initial set $\Lambda$ is in effect. Proof of sufficiency is complete.
\section{Equivalent problems}
It should be noted the technique of the proof of necessity in the theorem~\ref{T1} is valid provided the less restrictive assumptions are imposed. That allows to consider a plurality of equivalent relevant problems. Consider some exemplary situations.   

To begin with the multiplicities of nodes $\mu_k$ are irrelevant to the proof of necessity. As it is the properties (i) and (ii) are necessary for simple interpolation. On the other hand simple interpolation is a particular case of interpolation that takes into account multiplicities. 

Further let us consider the following problem of point-wise crude approximation:

{\noindent \it To acquire the conditions on $\mathcal M$ and $\Lambda$ that allow for some given sequence $\delta_k,\,k\in\bN,\,$ and for arbitrary data set $\{b_k,\,k\in\bN\},\,$ to find a function $u$ from $\Sigma (\bC,\Lambda)$ such that  $\left |u(\mu_k)-b_k\right |\leqslant \delta_k$ (or, putting this another way, for every entire function $f$ to find a function $u$ from $\Sigma (\bC,\Lambda)$ such that  $\left |u(\mu_k)-f(\mu_k)\right |\leqslant \delta_k$).}

The fact that interpolation by functions from $\Sigma (\bC,\Lambda)$  is operable (theorem~\ref{T1}) obviously implies solubility of the last problem. Conversely one readily notes that the slightly patched proof of necessity theorem~\ref{T1} gives the necessity of properties (i) and (ii) for solubility of the problem of point-wise crude approximation. 
In such a way we obtain the following.
\begin{theorem}\label{T3} The existence of representations~\ref{F1} (taking into account mutiplicities $m_k\geqslant 1$) is equivalent to both solubility of the problem of simple interpolation (all $m_k=1$) and solubility of the problem of point-wise crude approximation by functions from $\Sigma (\bC,\Lambda)$. So both problems are equivalent to properties (i) and (ii) of $\mathcal M$ and $\Lambda$ in the theorem~\ref{T1}.
\end{theorem}


\bigskip
\begin{thebibliography}{1}
\bibitem{Kri1} A. S. Krivosheev, “A fundamental principle for invariant subspaces in convex domains”, Izv. RAN. Ser. Mat., 68:2 (2004), 71–136
\bibitem{Kri2} A. S. Krivosheyev, O. A. Krivosheyeva, “A closedness of set of Dirichlet series sums”, Ufimsk. Mat. Zh., 5:3 (2013), 96–120
doi:10.13108/2013-5-3-94
%\bibitem{Krio1} {\it Кривошеева О.А. Область сходимости рядов экспоненциальных многочленов} // Уфимск. матем. журн. 5:4. 2013. C.~84–90. 
\bibitem{Kri4}	O.A. Krivosheeva,   A.S. Krivosheev,  Funct Anal Its Appl (2012) 46: 249. doi:10.1007/s10688-012-0033-1
\bibitem{Leonu} Leont'ev A.F. Series of exponentials, Moscow: Nauka. 1976. 536~p. (In Russian)
\bibitem{Leon} Leont'ev A.F. Sequences of polynomials of exponentials, Moscow: Nauka. 1980. 384~p. (In Russian)
\bibitem{Mcont} S.G. Merzlyakov, “Integrals of exponential functions with respect to Radon measure”, Ufa Math. Journal, 3:2 (2011), 56–78
\bibitem{MP} S.G. Merzlyakov, S.V. Popenov, “Interpolation with multiplicity by series of exponentials in H(C) with nodes on the real axis”, Ufa Math. Journal, 5:3 (2013), 127–140
\bibitem{MP1} S.G. Merzlyakov, S.V. Popenov, “Interpolation by series of exponentials in H(D) with real nodes”, Ufa Math. Journal, 7:1 (2015), 46–57
\bibitem{Mer1} S.G. Merzlyakov, “Invariant subspaces of the operator of multiple differentiation”, Math. Notes, 33:5 (1983), 361–368
\bibitem{Na1} V.V. Napalkov, “Complex Analysis and the Cauchy Problem for Convolution Operators”, Proc. Steklov Inst. Math., 235 (2001), 158–161
\bibitem{Na} V. V. Napalkov, A.A. Nuyatov, “The multipoint de la Valle\'e Poussin problem for a convolution operator”, Sb. Math., 203:2 (2012), 224–233
\bibitem{NN} V. V. Napalkov, A.A. Nuyatov, “Multipoint Valle\'e Poussin problem for convolution operators with nodes defined inside an angle”, Theoret. and Math. Phys., 180:2 (2014), 983–989 
\bibitem{NZ} K.R. Zimens, V.V. Napalkov, “The multiple de la Valle\'e Poussin problem on convex domains in the kernel of the convolution operator”, Dokl. Math., 90:2 (2014), 581-583 doi:10.1134/S1064562414060234
\bibitem{NapP} V.V. Napalkov, S.V. Popenov, “The holomorphic Cauchy problem for convolution operators in analytically uniform spaces and Fischer's representations”,
Doklady. Mathematics, 2001. 
 \bibitem{Seb} Sebasti\~ao e Silva J. “Su serte classi di spazi localmente convessi importanti per le applicazioni“. Rendiconti de mathematica e delle sue applicazoni. Roma. (5) 14 (1955),388-410.
 \bibitem{Rud} Rudin W. Functional analysis, New York: McGRAW-HILL COMPANY. 1973..
\bibitem{10} H.~Muggli. “Differentialgleichungen unendlich
hoher Ordnung mit konstanten Koeffizienten”. //  Comment.  Math. Helv. V. 11 (1938), 151-179. 
\bibitem{Her}  H\"ormander L. An introduction to complex analysis of several variables. Princenton, New Jersey: D. VAN NOSTRAND COMPANY, INC. 1966.
\end{thebibliography}
\end{document}